\documentclass[12pt]{amsart}
\usepackage{amssymb,hyperref}

\newtheorem{theorem}{Theorem}[section]
\newtheorem{definition}[theorem]{Definition}

\newtheorem{proposition}[theorem]{Proposition}

\newtheorem{lemma}[theorem]{Lemma}

\begin{document}

\title[Twisted group algebras and free hypergeometric laws]{Quantum automorphisms of twisted group algebras and free hypergeometric laws}

\author{Teodor Banica}
\address{T.B.: Department of Mathematics, Cergy-Pontoise University, 95000 Cergy-Pontoise, France. {\tt teodor.banica@u-cergy.fr}}

\author{Julien Bichon}
\address{J.B.: Department of Mathematics, Clermont-Ferrand University, Campus des Cezeaux, 63177 Aubiere Cedex, France. {\tt bichon@math.univ-bpclermont.fr}}

\author{Stephen Curran}
\address{S.C.: Department of Mathematics, University of California, Berkeley, CA 94720, USA. {\tt curransr@math.berkeley.edu}}


\begin{abstract}
We prove that we have an isomorphism of type $A_{aut}(\mathbb C_\sigma[G])\simeq A_{aut}(\mathbb C[G])^\sigma$, for any finite group $G$, and any $2$-cocycle $\sigma$ on $G$. In the particular case $G=\mathbb Z_n^2$, this leads to a Haar-measure preserving identification between the subalgebra of $A_o(n)$ generated by the variables $u_{ij}^2$, and the subalgebra of $A_s(n^2)$ generated by the variables $X_{ij}=\sum_{a,b=1}^np_{ia,jb}$.  Since $u_{ij}$ is ``free hyperspherical'' and $X_{ij}$ is ``free hypergeometric'', we obtain in this way a new free probability formula, which at $n=\infty$ corresponds to the well-known relation between the semicircle law, and the free Poisson law.
\end{abstract}

\maketitle

\section*{Introduction}

The notion of quantum automorphism group was introduced about 10 years ago, in Wang's paper \cite{wan}. The idea is as follows:

\begin{enumerate}
\item First, a finite quantum space $X$ is by definition the ``spectrum'' of a finite dimensional    $\mathbb C^*$-algebra $A$. This space comes with the measure corresponding to the ``canonical trace'' $tr:A\to\mathbb C$, obtained via the left regular representation.

\item Wang proved that this space $X$ has a quantum automorphism group $G_{aut}(X)$. In algebraic terms, the result is that there is a universal Hopf algebra $A_{aut}(A)$, coacting on $A$, and leaving the canonical trace invariant.
\end{enumerate}

As a basic, motivating example, for the space $X=\{n\ {\rm points}\}$ we have $A=\mathbb C^n$. Just by using some examples of dual coactions, one can prove then that for $n\geq 4$ the algebra $A_{aut}(A)$ is not commutative, and infinite dimensional. In other words, the space consisting of $n\geq 4$ points has an infinite number of ``quantum permutations''. See Wang \cite{wan}.

The quantum automorphism groups, and their quantum subgroups, were systematically investigated in the last years. See e.g. \cite{ban}, \cite{bbi}, \cite{bgs}, \cite{bi1}, \cite{bi2}.

In this paper we establish a general isomorphism result, of the following type:
$$A_{aut}(\mathbb C_\sigma[G])\simeq A_{aut}(\mathbb C[G])^\sigma$$

Here $G$ is a finite group, $\mathbb C[G]$ is its group algebra, $\sigma$ is a 2-cocycle on $G$, and $\mathbb C_\sigma[G]$ is the corresponding twisted group algebra. As for the algebra on the right, this is a 2-cocycle twist. In other words, our result is that ``the quantum automorphism group of a twisted group algebra is the twist of the original quantum automorphism group''.

As an illustration for this result, we will work out in detail the case of the abelian group $G=\mathbb Z_n^2$, coming with the cocycle $\sigma((i,j),(k,l))=w^{jk}$, where $w=e^{2\pi i/n}$. Here $A_{aut}(\mathbb C[G])$ is Wang's quantum permutation algebra $A_s(n^2)=A_{aut}(\mathbb C^{n^2})$, and $A_{aut}(\mathbb C_\sigma[G])$ is Wang's projective quantum rotation algebra $PA_o(n)$. Thus, the above isomorphism reads: 
$$PA_o(n)\simeq A_s(n^2)^\sigma$$

In other words, the above general algebraic result leads in this case to a non-trivial relation between quantum rotations and quantum permutations!

We will further investigate this phenomenon, by using concepts from Voiculescu's free probability theory \cite{vdn}. The idea is that, probabilistically speaking, the above result tells us that the variables $u_{ij}^2\in A_o(n)$ have the same law as the following variables: 
$$X_{ij}=\frac{1}{n}\sum_{a,b=1}^np_{ia,jb}\in A_s(n^2)$$

Here $u_{ij}$ and $p_{ia,jb}$ are respectively the standard generators of $A_o(n)$ and $A_s(n^2)$.

The point now is that the variables $u_{ij}$ can be regarded as being ``free hyperspherical variables'', and the variables $X_{ij}$ can be regarded as being (rescaled) ``free hypergeometric variables''. So, what we have here is a new free probability formula, which at $n=\infty$ corresponds to the well-known relation between the semicircle law, and the free Poisson law. We will make several comments on this result, including a further investigation of the free hypergeometric laws, for more general values of the parameters.

The paper is organized as follows: in 1 we state and prove the main algebraic result, in 2 we work out the case $G=\mathbb Z_n^2$, and in 3 we discuss the free probability aspects. The final section, 4, contains some remarks, comments, and open questions.

\subsection*{Acknowledgements}

It is a pleasure to thank B. Collins and R. Speicher for several useful discussions, prior to the work leading to the present article. The work of T.B. and J.B. was supported by the ANR grant ``Galoisint''.

\section{The twisting result} 

Let $A$ be a finite dimensional $\mathbb C^*$-algebra, i.e. a finite direct sum of matrix algebras. We recall from \cite{ban} that the canonical trace $tr:A\to\mathbb C$ is obtained by composing the left regular representation $A\to\mathcal L(A)$ with the usual (normalized) trace $\mathcal L(A)\to\mathbb C$. 

Note that in the cases $A=\mathbb C^n$ and $A=M_n(\mathbb C)$ we obtain the usual trace. 

We let $A_{aut}(A)$ be the universal Hopf algebra coacting on the algebra $A$, such that the canonical trace $tr$ is invariant. This construction is the $\phi=tr$ particular case of Wang's general construction $A_{aut}(A,\phi)$, from \cite{wan}. Regarding the choice $\phi=tr$, see \cite{ban}.

In this paper we are interested in the case where $A=\mathbb C[G]$ is the convolution algebra of a finite group $G$. If $\{e_g\}$ denotes the standard basis of $\mathbb C[G]$, we have $tr(e_g)=\delta_{g1}$.

\begin{proposition}
$A_{aut}(\mathbb C[G])$ is isomorphic to the abstract algebra presented by generators $x_{g,h}$ with $g,h\in G$, with the following relations:
$$x_{1,g}=x_{g,1}=\delta_{1g},\quad x_{s,gh}=\sum_{t\in G}x_{st^{-1},g}x_{t,h},\quad 
x_{gh,s}=\sum_{t\in G}x_{g,t^{-1}}x_{h,ts}$$
The comultiplication, counit and antipode are given by the following formulae:
$$\Delta(x_{g,h})=\sum_{s\in G}x_{g,s}\otimes x_{s,h},\quad\varepsilon(x_{g,h})=\delta_{gh},\quad S(x_{g,h})=x_{h^{-1},g^{-1}}$$
\end{proposition}

\begin{proof}
This follows from a direct verification: in fact, for $A=\mathbb C[G]$, the algebra $A_{aut}(A)$ constructed by Wang in \cite{wan} is precisely the universal algebra in the statement.
\end{proof}

Now let $\sigma:G\times G\to\mathbb C^*$ be a normalized 2-cocycle on $G$, i.e. a map satisfying $\sigma(gh,s)\sigma(g,h)=\sigma(g,hs)\sigma(h,s)$ and $\sigma(g,1)=\sigma(1,g)=1$. The twisted group algebra $\mathbb C_\sigma[G]$ is by definition the vector space $\mathbb C[G]$, with product $e_ge_h=\sigma(g,h)e_{gh}$. 

It is routine to check that the canonical trace of $\mathbb C_\sigma[G]$ is the same as the one of ${\mathbb C}[G]$. We have the following generalization of Proposition 1.1.

\begin{proposition}
$A_{aut}(\mathbb C_\sigma[G])$ is isomorphic to the abstract algebra presented by generators $x_{g,h}$ with $g,h\in G$, with the relations $x_{1,g}=x_{g,1}=\delta_{1g}$ and:
$$\sigma(g,h)x_{s,gh}=\sum_{t\in G}\sigma(st^{-1},t)x_{st^{-1},g}x_{t,h},\quad \sigma(g,h)^{-1}x_{gh,s}=\sum_{t\in G}\sigma(t^{-1},ts)^{-1}x_{g,t^{-1}}x_{h,ts}$$
The comultiplication, counit and antipode are given by the following formulae:
$$\Delta(x_{g,h})=\sum_{s\in G}x_{g,s}\otimes x_{s,h},\quad \varepsilon(x_{g,h})=\delta_{gh},\quad S(x_{g,h})=\sigma(h^{-1},h)\sigma(g^{-1},g)^{-1}x_{h^{-1},g^{-1}}$$
\end{proposition}

\begin{proof}
Once again, this follows from a direct verification. Note that by using cocycle identities we obtain $\sigma(g,g^{-1})=\sigma(g^{-1},g)$, needed in the proof.
\end{proof}

In what follows, we will prove that the Hopf algebras $A_{aut}({\mathbb C}[G])$ and $A_{aut}({\mathbb C}_\sigma[G])$ are related by a ``cocycle twisting'' operation. Let us begin with some preliminaries.

Let $H$ be a Hopf algebra. We use the Sweedler notation $\Delta(x)=\sum x_1\otimes x_2$. Recall (see e.g. \cite{doi}) that a left 2-cocycle is a convolution invertible linear map
$\sigma:H\otimes H\to\mathbb C$ satisfying $\sigma(x_1,y_1)\sigma(x_2y_2,z)=\sigma(y_1,z_1)\sigma(x,y_2z_2)$ and $\sigma(x,1)=\sigma(1,x)=\varepsilon(x)$, for any $x,y,z \in H$. Note that $\sigma$ is a left 2-cocycle if and only if $\sigma^{-1}$, the convolution inverse of $\sigma$, is a right 2-cocycle, in the sense that we have $\sigma^{-1}(x_1y_1, z)\sigma^{-1}(x_1,y_2)=\sigma^{-1}(x,y_1z_1)\sigma^{-1}(y_2,z_2)$ and $\sigma^{-1}(x,1)=\varepsilon(x)=\sigma^{-1 }(1)$, for any $x,y,z \in H$.

Given a left 2-cocycle $\sigma$ on $H$, one can form the 2-cocycle twist $H^\sigma$ as follows. As a coalgebra, $H^\sigma=H$, and an element $x\in H$, when considered in $H^\sigma$, is denoted $[x]$. The product in $H^\sigma$ is defined, in Sweedler notation, by: 
$$[x][y]=\sum\sigma(x_1,y_1)\sigma^{-1}(x_3,y_3)[x_2y_2]$$

Note that the cocycle condition ensures the fact that we have indeed a Hopf algebra. For the formula of the antipode, that we will not need here, see \cite{doi}.

Note that the coalgebra isomorphism $H\to H^\sigma$ given by $x\to [x]$ commutes with the respective Haar integrals (as soon as $H$ has a Haar integral, of course).

We are now in position to state and prove our main theorem.

\begin{theorem}
If $G$ is a finite group and $\sigma$ is a $2$-cocycle on $G$, the Hopf algebras $A_{aut}(\mathbb C[G])$ and $A_{aut}(\mathbb C_\sigma[G])$ are $2$-cocycle twists of each other.
\end{theorem}

\begin{proof}
We use the Hopf algebra map $\pi:A_{aut}(\mathbb C[G])\to\mathbb C[G]$, given by $x_{g,h}\to\delta_{gh}e_g$.

Our 2-cocycle $\sigma:G\times G\to\mathbb C^*$ can be extended by linearity into a linear map $\sigma:\mathbb C[G]\otimes\mathbb C[G]\to\mathbb C$, which, since $\mathbb C[G]$ is a group algebra, is a left and right 2-cocycle in the above sense. Consider the following map:
$$\alpha=\sigma(\pi \otimes \pi):A_{aut}(\mathbb C[G])\otimes A_{aut}(\mathbb C[G]) \to\mathbb C[G]\otimes\mathbb C[G]\to\mathbb C$$

Then $\alpha$ is a left and right 2-cocycle, because it is induced by a cocycle on a group algebra, and so is its convolution inverse $\alpha^{-1}$. Thus we can construct the twisted algebra $A_{aut}({\mathbb C}[G])^{\alpha^{-1}}$, and in this algebra we have:
$$[x_{g,h}][x_{r,s}]=\alpha^{-1}(x_g,x_r)\alpha(x_h,x_s)[x_{g,h}x_{r,s}]=\sigma(g,r)^{-1}\sigma(h,s)[x_{g,h}x_{r,s}]$$

By using this, we see that:
\begin{eqnarray*}
\sum_{t \in G}\sigma(st^{-1},t)[x_{st^{-1},g}][x_{t,h}] 
&=&\sum_{t \in G}\sigma(st^{-1},t)\sigma(st^{-1},t)^{-1}\sigma(g,h)[x_{st^{-1},g}x_{t,h}] \\
&=&\sigma(g,h)[x_{s,gh}]
\end{eqnarray*}

Similarly, we have:
$$\sum_{t\in G}\sigma(t^{-1},ts)^{-1}[x_{g,t^{-1}}][x_{h,ts}]=\sigma(g,h)^{-1}[x_{gh,s}]$$

Thus there exists a Hopf algebra map $F:A_{aut}(\mathbb C_\sigma[G])\to A_{aut}(\mathbb C[G])^{\alpha^{-1}}$, given by $x_{g,h}\to [x_{g,h}]$. This map is clearly surjective, and is injective by a standard fusion semiring argument (both Hopf algebras have the same fusion semiring by \cite{ban}).
\end{proof}

Associated with any $2$-cocycle are the following quantities:
$$\Omega(g_1,\ldots,g_m)=\prod_{k=1}^{m-1}\sigma(g_1\ldots g_k,g_{k+1})$$

With this notation, we have the following technical reformulation of Theorem 1.3. 

\begin{proposition}
If $G$ is a finite group and $\sigma$ is a $2$-cocycle on $G$, then
$$F(x_{g_1,h_1}\ldots x_{g_m,h_m})=\Omega(g_1,\ldots,g_m)^{-1}\Omega(h_1,\ldots,h_m)x_{g_1,h_1}\ldots x_{g_m,h_m}$$
is a coalgebra isomorphism $A_{aut}(\mathbb C_\sigma[G])\simeq A_{aut}(\mathbb C[G])$, commuting with the Haar integrals.
\end{proposition}

\begin{proof}
This is indeed just a technical reformulation of Theorem 1.3.
\end{proof}

\begin{theorem}
Let $X\subset G$ be such that $\sigma(g,h)=1$ for any $g,h \in X$, and consider the subalgebra $B_X\subset A_{aut}(\mathbb C_\sigma[G])$ generated by the elements $x_{g,h}$, with $g,h\in X$. Then we have an injective algebra map $F_0:B_X\to A_{aut}(\mathbb C[G])$, given by $x_{g,h}\to x_{g,h}$.
\end{theorem}

\begin{proof}
With the notations in the proof of Theorem 1.3, we have the following equality in $A_{aut}(\mathbb C[G])^{\alpha^{-1}}$, for any $g_i,h_i,r_i,s_i\in X$:
$$[x_{g_1,h_1}\ldots x_{g_p,h_p}] \cdot [x_{r_1,s_1}\ldots x_{r_q,s_q}] 
= [x_{g_1,h_1}\ldots x_{g_p,h_p}x_{r_1,s_1}\ldots x_{r_q,s_q}]$$

Now $F_0$ can be defined to be the composition of $F_{|B_X}$ with the linear isomorphism $A_{aut}(\mathbb C[G])^{\alpha^{-1}}\to A_{aut}(\mathbb C[G])$, $[x]\to x$, and is clearly an injective algebra map.   
\end{proof}

\section{Rotations and permutations}

In this section we discuss some concrete consequences of the general results established in the previous section. These will concern Wang's quantum permutation groups \cite{wan}.

Consider the additive group $G= \mathbb Z_n^2$. Let  $w\in\mathbb C^*$ be a primitive $n$-th root of unity, and consider the map $\sigma:\mathbb Z_n^2\times\mathbb Z_n^2\to\mathbb C^*$ given by:
$$\sigma((i,j),(k,l))=w^{jk}$$ 

It is easy to see that $\sigma$ is a bicharacter, and hence a 2-cocycle on $\mathbb Z_n^2$.

We denote by $E_{ij}$ with $i,j \in\mathbb Z_n$ the standard elementary matrices in $M_n(\mathbb C)$.

\begin{lemma}
The linear map given by
$$\psi(e_{(i,j)})=\sum_{k=0}^{n-1}{w}^{ki}E_{k,k+j}$$
defines an isomorphism of algebras $\psi:\mathbb C_\sigma[\mathbb Z_n^2]\simeq M_n(\mathbb C)$. 
\end{lemma}

\begin{proof}
Consider indeed the following linear map:
$$\psi'(E_{ij})=\frac{1}{n}\sum_{k=0}^{n-1}{w}^{-ik}e_{(k,j-i)}$$
 
It is routine to check that both $\psi,\psi'$ are morphisms of algebras, and that these maps are inverse to each other. In particular, $\psi$ is an isomorphism of algebras, as stated.
\end{proof}

\begin{lemma}
The algebra map given by
$$\varphi(u_{ij}u_{kl}) = \frac{1}{n}\sum_{a,b=0}^{n-1}{w}^{ai-bj}x_{(a,k-i),(b,l-j)}$$
defines a Hopf algebra isomorphism $\varphi:A_{aut}(M_n(\mathbb C))\simeq A_{aut}(\mathbb C_\sigma[\mathbb Z_n^2])$.
\end{lemma}

\begin{proof}
Consider the universal coactions on the two algebras in the statement:
\begin{eqnarray*}
 \alpha:M_n(\mathbb C)&\to&M_n({\mathbb C})\otimes A_{aut}(M_n({\mathbb C}))\\
 \beta : {\mathbb C}_\sigma[{\mathbb Z}_n^2]&\to&{\mathbb C}_\sigma[{\mathbb Z}_n^2] \otimes A_{aut}({\mathbb C}_\sigma[{\mathbb Z}_n^2])
 \end{eqnarray*}
 
 In terms of the standard bases, these coactions are given by:
 \begin{eqnarray*}
 \alpha(E_{ij})&=&\sum_{kl}E_{kl}\otimes u_{ki}u_{lj}\\
 \beta(e_{(i,j)})&=&\sum_{kl} e_{(k,l)}\otimes x_{(k,l),(i,j)}
 \end{eqnarray*}

 We use now the identification $\mathbb C_\sigma[\mathbb Z_n^2]\simeq M_n(\mathbb C)$ given by Lemma 2.1. The resulting coaction $\gamma:M_n(\mathbb C)\to M_n(\mathbb C)\otimes A_{aut}(\mathbb C_\sigma[\mathbb Z_n^2])$ is then given by the following formula:
 $$\gamma(E_{ij})=\frac{1}{n}\sum_{ab}E_{ab}\otimes\sum_{kr}w^{ar-ik} x_{(r,b-a),(k,j-i)}$$

By comparing with the formula of $\alpha$, we obtain the isomorphism in the statement.
\end{proof}

\begin{lemma}
The algebra map given by
$$\rho(x_{(a,b),(i,j)})=\frac{1}{n^2}\sum_{klrs}w^{ki+lj-ra-sb}p_{(r,s),(k,l)}$$
defines a Hopf algebra isomorphism $\rho:A_{aut}(\mathbb C[\mathbb Z_n^2])\simeq
 A_{aut}(C(\mathbb Z_n^2))$.
\end{lemma}

\begin{proof}
This is similar to the proof of the previous lemma, by using the Fourier transform isomorphism $\mathbb C[\mathbb Z_n^2]\simeq C(\mathbb Z_n^2)$.
\end{proof}

Consider now the Wang algebras $A_o(n)$ and $A_s(n^2)$, with standard generators denoted $(u_{ij})_{i,j=1,\ldots,n}$ and $(p_{ia,jb})_{i,j,a,b=1,\ldots,n}$. That is, $u=(u_{ij})$ is the universal $n\times n$ orthogonal matrix, and $p=(p_{ia,jb})$ is the universal $n^2\times n^2$ magic unitary matrix. See \cite{wan}. 

We recall that we have canonical identifications, as follows:
\begin{eqnarray*}
A_s(n^2)&=&A_{aut}(\mathbb C^{n^2})\\
PA_o(n)&=&A_{aut}(M_n(\mathbb C))
\end{eqnarray*}

Here the projective version of a pair $(A,u)$ is by definition the pair $(PA,v)$, where $v=u\otimes\bar{u}$ and $PA=<v_{ij}>$. For full details regarding the above equalities, see \cite{ban}.

\begin{theorem}
Let $n\geq 2$ and $w\in\mathbb C^*$ be a primitive $n$-th root of unity. Then
$$\Theta(u_{ij}u_{kl})=\frac{1}{n}\sum_{ab=0}^{n-1}w^{-a(k-i)+b(l-j)}p_{ia,jb}$$
defines a coalgebra isomorphism $PA_o(n)\to A_s(n^2)$, commuting with the Haar integrals.
 \end{theorem}
 
 \begin{proof}
This follows from the general isomorphism results in Theorem 1.3 and Proposition 1.4, by combining them with the various isomorphisms from the lemmas above.
\end{proof}

\begin{theorem}
The following two algebras are isomorphic, via $u_{ij}^2\to X_{ij}$:
\begin{enumerate}
\item The algebra generated by the variables $u_{ij}^2\in A_o(n)$.

\item The algebra generated by $X_{ij}=\frac{1}{n}\sum_{a,b=1}^np_{ia,jb}\in A_s(n^2)$
\end{enumerate}
\end{theorem}

\begin{proof}
We have $\Theta(u_{ij}^2)=X_{ij}$, so it remains to prove that if $B$ is the subalgebra of $A_{aut}(M_n(\mathbb C))$ generated by the variables $u_{ij}^2$, then $\Theta_{|B}$ is an algebra morphism.

We let $X=\{(i,0)|i\in\mathbb Z_n\}\subset\mathbb Z_n^2$. Then $X$ satisfies to the assumption in Theorem 1.5, and $\varphi(B) \subset B_X$. Thus by Theorem 1.5, $\Theta_{|B}=\rho F_0\varphi_{|B}$ is indeed an algebra morphism.
\end{proof}

\section{Free hypergeometric laws}

Let $(A,\varphi)$ be a complex algebra, coming with a linear form $\varphi:A\to\mathbb C$. We recall that the law (or distribution) of an element $a\in A$ is the linear form $\mu_a:\mathbb C[X]\to\mathbb C$ given by $P\to\varphi(P(a))$. More generally, the joint law (or distribution) of a family $a=(a_1,\ldots,a_n)$ is the linear form $\mu_a:\mathbb C<X_1,\ldots,X_n>\to\mathbb C$ given by $P\to\varphi(P(a))$. Thus, saying that two families $a_1,\ldots,a_n\in A$ and $a_1',\ldots,a_n'\in A'$ have the same law is the same as saying that we have the following equalities, for any noncommutative polynomial $P$:
$$\varphi(P(a_1,\ldots,a_n))=\varphi'(P(a_1',\ldots,a_n'))$$

With this notation, we know from Theorem 2.5 that the variables $u_{ij}^2\in A_o(n)$ have the same joint law as the variables $X_{ij}\in A_s(n^2)$. In this section we will put this result into a more general framework, related to Voiculescu's free probability theory \cite{vdn}.

Let us begin by giving an independent proof for the above equality of distributions. We use the Weingarten formula \cite{bco}. We recall that the representation-theoretic sets of partitions for the algebras $A=A_o(n),A_s(n^2)$ are respectively the set of noncrossing pairings $NC_2(2k)$, and the set of noncrossing partitions $NC(k)$. Associated to each of these sets are the Gram matrix $G_m(\pi,\sigma)=m^{|\pi\vee\sigma|}$ and the Weingarten matrix $W_m=G_m^{-1}$, where $\vee$ is the join operation, and $m=n,n^2$ respectively. The Haar functional of $A$ can be computed explicitely in terms of $W_m$. For full details here, see \cite{bcs}.

\begin{lemma}
The Weingarten matrices of $A_o(n)$ and $A_s(n^2)$ are related by
$$W_{NC_2(2k),n}(\pi,\sigma)=n^{|\tilde{\pi}|+|\tilde{\sigma}|-k}W_{NC(k),n^2}(\tilde{\pi},\tilde{\sigma})$$
where $\pi\to\tilde{\pi}$ is the ``cabling'' operation, obtained by collapsing neighbors.
\end{lemma}

\begin{proof}
We use the following general formula, due to Kodiyalam and Sunder \cite{ksu}:
$$|\pi\vee\sigma|=k+2|\tilde{\pi}\vee\tilde{\sigma}|-|\tilde{\pi}|-|\tilde{\sigma}|$$

See also \cite{csp}. Now in terms of Gram matrices, we obtain:
$$G_{NC_2(2k),n}=n^{k-|\tilde{\pi}|-|\tilde{\sigma}|}G_{NC(k),n^2}(\tilde{\pi},\tilde{\sigma})$$

By taking the inverse, this gives the formula in the statement.
\end{proof}

\begin{theorem}
The family of variables
$$\left\{X_{ij}=\frac{1}{n}\sum_{a,b=1}^np_{ia,jb}\right\}\subset A_s(n^2)$$
has the same law as the family of variables $\{u_{ij}^2\}\subset A_o(n)$.
\end{theorem}

\begin{proof}
We use the Weingarten formula, which gives:
\begin{eqnarray*}
\int X_{ij}^k
&=&\sum_{\pi,\sigma\in NC_2(2k)}n^{|\tilde{\pi}|+|\tilde{\sigma}|-k}W_{NC(k),n^2}(\tilde{\pi},\tilde{\sigma})\\
\int u_{ij}^{2k}
&=&\sum_{\pi,\sigma\in NC_2(2k)}W_{NC_2(2k),n}(\pi,\sigma)
\end{eqnarray*}

By Lemma 3.1 the terms on the right are equal, and this gives the equality of laws for the individual variables. For the general statement, the proof is similar.
\end{proof}

The variables $X_{ij}$ appearing in Theorem 3.2 have the following generalization.

\begin{definition}
The noncommutative random variable
$$X(n,m,N) = \sum_{i=1}^n \sum_{j=1}^m u_{ij} \in A_s(N)$$
is called free hypergeometric, of parameters $(n,m,N)$.
\end{definition}

The terminology here comes from the fact that the variable $X'(n,m,N)$, defined as above, but over the algebra $C(S_n)$, follows a hypergeometric law of parameters $(n,m,N)$.

In general, the free hypergeometric laws seem to be quite difficult to compute. A first result in this direction, heavily relying on a result recently obtained in \cite{bcz}, is as follows.

\begin{theorem}
The moments of $X(n,n,n^2)$ are given by
$$\int X(n,n,n^2)^k\,dx=\frac{n^k}{(n+1)^k}\cdot\frac{q+1}{q-1}\cdot\frac{1}{k+1}\sum_{r=-k-1}^{k+1}(-1)^r\begin{pmatrix}2k+2\cr k+r+1\end{pmatrix}\frac{r}{1+q^r}$$
where $q\in [-1,0)$ is given by $q+q^{-1}=-n$.
\end{theorem}

\begin{proof}
First, $X(n,n,n^2)/n$ is the variable $X_{ij}$ appearing by Theorem 3.2, having the same law as the variable $u_{ij}^2\in A_o(n)$. Now it is known from \cite{bdv} that $A_o(n)$ is monoidally equivalent to $C(SU_q(2))$, and, as explained in \cite{bcz}, one can use this fact for modelling $u_{ij}\in A_o(n)$ by a certain variable over $SU_q(2)$. This latter variable can be studied by using advanced calculus methods, and this leads to the above formula. See \cite{bcz}.
\end{proof}

As a first observation, the above result, or rather a version of it, namely Theorem 5.3 in \cite{bcz}, shows that the variables $X(n,n,n^2)$ superconverge with $n\to\infty$. For more about the superconvergence phenomenon in free probability, see Bercovici and Voiculescu \cite{bvo}.

The second result, that we would like to present now, is an exploration of the basic  asymptotic properties of the free hypergeometric laws.

\begin{theorem}
The free hypergeometric laws have the following properties:
\begin{enumerate}
\item Let $n,m,N\to\infty$, with $\frac{nm}{N}\to\lambda\in (0,\infty)$.  Then the law of $X(n,m,N)$ converges to the free Poisson law of parameter $\lambda$.

\item Let $n,m,N\to\infty$, with $\frac{n}{N}\to\nu\in (0,1)$ and $\frac{m}{N}\to 0$. Then the law of $S(n,m,N)=(X(n,m,N) - m\nu)/\sqrt{m\nu(1-\nu)}$ converges to a $(0,1)$-semicircle law.
\end{enumerate}
\end{theorem}

\begin{proof}
(1) From the Weingarten formula, we have:
$$\int X(n,m,N)^p = \sum_{\pi,\sigma \in NC(p)}W_{NC(p),N}(\pi,\sigma) n^{|\pi|}m^{|\sigma|}$$

Now, as explained for instance in \cite{cur}, we have:
$$W_{NC(p),N}(\pi,\sigma) = \begin{cases} N^{-|\pi|} + O(N^{-|\pi|-1}), & \pi = \sigma \\ O(N^{|\pi \vee \sigma| - |\pi| - |\sigma|}), & \pi \neq \sigma \end{cases}$$

It follows that we have:
$$W_{NC(p),N}(\pi,\sigma) n^{|\pi|}m^{|\sigma|} \to \begin{cases} \lambda^{|\pi|}, & \pi = \sigma\\ 0, & \pi \neq \sigma \end{cases}$$

Thus the $p$-th moment of $X(n,m,N)$ converges to $\sum_{\pi \in NC(p)} \lambda^{|\pi|}$, which is the $p$-th moment of the free Poisson distribution with parameter $\lambda$, and we are done.

(2) We need to show that the free cumulants satisfy:
$$\kappa^{(p)}[S(n,m,N),\ldots,S(n,m,N)]  \to \begin{cases} 1, & p = 2 \\ 0, & p \neq 2 \end{cases}$$

The case $p =1$ is trivial, so suppose $p\geq 2$.  We have:
$$\kappa^{(p)}[S(n,m,N),\ldots,S(n,m,N)]=(m\nu(1-\nu))^{-p/2}\kappa^{(p)}[X(n,m,N),\ldots,X(n,m,N)]$$

On the other hand, from the Weingarten formula, we have:
\begin{eqnarray*}
&&\kappa^{(p)}[X(n,m,N),\dotsc,X(n,m,N)] \\
&=&\sum_{{w} \in NC(p)} \mu_p({w},1_p) \prod_{V \in {w}} \sum_{\pi_V,\sigma_V \in NC(V)} W_{NC(V),N}(\pi_V,\sigma_V)n^{|\pi_V|}m^{|\sigma_V|}\\
&=&\sum_{{w} \in NC(p)} \mu_p({w},1_p) \prod_{V \in {w}} \sum_{\pi_V,\sigma_V \in NC(V)} (N^{-|\pi_V|}\mu_{|V|}(\pi_V,\sigma_V) + O(N^{-|\pi_V|-1}))n^{|\pi_V|}m^{|\sigma_V|}\\
&=&\sum_{\substack{\pi,\sigma \in NC(p)\\ \pi \leq \sigma}} (N^{-|\pi|}\mu_p(\pi,\sigma) + O(N^{-|\pi|-1}))n^{|\pi|}m^{|\sigma|} \sum_{\substack{{w} \in NC(p)\\ \sigma \leq {w}}} \mu_p({w},1_p)
\end{eqnarray*}

We use now the following standard identity: 
$$\sum_{\substack{{w} \in NC(p)\\ \sigma \leq {w}}} \mu_p({w},1_p) = \begin{cases} 1, & \sigma = 1_p\\ 0, & \sigma \neq 1_p \end{cases}$$

This gives the following formula for the cumulants:
\begin{equation*}
  \kappa^{(p)}[X(n,m,N),\ldots,X(n,m,N)]=m\sum_{\pi \in NC(p)} (N^{-|\pi|}\mu_p(\pi,1_p) + O(N^{-|\pi|-1}))n^{|\pi|}
\end{equation*}

It follows that for $p \geq 3$ we have, as desired:
$$\kappa^{(p)}[S(n,m,N),\dotsc,S(n,m,N)]\to 0$$

As for the remaining case $p = 2$, here we have:
\begin{eqnarray*}
 \kappa^{(2)}[S(n,m,N),S(n,m,N)] 
&\to&\frac{1}{\nu(1-\nu)} \sum_{\pi \in NC(2)} \nu^{|\pi|}\mu_2(\pi,1_2)\\
&=&\frac{1}{\nu(1-\nu)}\bigl(\nu - \nu^2)\\
&=&1
\end{eqnarray*}

This gives the result.
\end{proof}

\section{Concluding remarks}

We have seen in this paper that Wang's quantum automorphism groups in \cite{wan} are subject to some general twisting results, and that these results are of relevance in the general context of Voiculescu's free probability theory \cite{vdn}. Several questions appear: 
\begin{enumerate}
\item What is the most general twisting result for quantum automorphism groups? An answer here is $A_{aut}(H_\sigma)\simeq A_{aut}(H)^\sigma$, for any finite dimensional Hopf algebra $H$. This result, whose proof is much more technical, will be discussed somewhere else.

\item Besides the quantum rotation/permutation application, is there any other concrete, probabilistic application of our general twisting result? One interesting example here seems to be $G=(\mathbb Z_2\times\mathbb Z_2)^n$, leading to the algebra $M_{2^n}(\mathbb C)$.

\item Can one include the monoidal equivalence used in \cite{bcz} into the above considerations? The point is that the variable $u_{ij}\in A_o(n)$ has the same law as a certain variable over $SU_q(2)$, where $q+q^{-1}=-n$, so the result in \cite{bcz} can be probably stated and proved by using $A_s(n^2)$ and $SU_q(2)$ only. However, it is not clear how to do so.

\item Are there any other free hypergeometric laws, that can be explicitely computed by using $A_o(F)$? In principle the answer here is no, first because of the concluding remarks in \cite{bcz}, and second, because of the ``no-atoms'' results of Voigt in \cite{voi}.

\item Do we have superconvergence, in the sense of Bercovici and Voiculescu \cite{bvo}, to the limiting distributions in Theorem 3.5? Note that in the case $n=m=\sqrt{N}$ the superconvergence appears indeed, as explained after Theorem 3.4.
\end{enumerate}

Finally, let us point out the fact that, as explained in \cite{bgo}, the free hyperspherical laws appear in connection with the study of a number of interesting ``noncommutative spaces'', such as the free and half-liberated spheres, or projective spaces. We do not know yet if the relation with the free hypergeometric laws, that we found in this paper, can be of help here, but we intend to come back to this question in some future work.

\end{document}